\documentclass[a4paper,12pt]{amsart}
\usepackage{amsfonts}
\usepackage{amsmath,amssymb,amsthm}

\everymath{\displaystyle}

\newcommand{\E}{\mathbb E}
\newcommand{\R}{\mathbb R}
\newcommand{\tr}{\mathrm{tr}}
\newcommand{\ds}{\displaystyle}

\newtheorem{thm}{Theorem}[section]
\newtheorem{cor}[thm]{Corollary}

\newtheorem{prop}[thm]{Proposition}
\theoremstyle{definition}

\theoremstyle{remark}
\newtheorem{rem}[thm]{Remark}

\begin{document}

\title[General Rotational Surfaces in Pseudo-Euclidean 4-Space]
{General Rotational Surfaces in Pseudo-Euclidean 4-Space with Neutral Metric}

\author{Yana  Aleksieva, Velichka Milousheva, Nurettin Cenk Turgay}

\address{Faculty of Mathematics and Informatics, Sofia University,
5 James Bourchier blvd., 1164 Sofia, Bulgaria}
\email{yana\_a\_n@fmi.uni-sofia.bg}

\address{Institute of Mathematics and Informatics, Bulgarian Academy of Sciences,
Acad. G. Bonchev Str. bl. 8, 1113, Sofia, Bulgaria;   "L.
Karavelov" Civil Engineering Higher School, 175 Suhodolska Str., 1373 Sofia, Bulgaria}
\email{vmil@math.bas.bg}

\address{Istanbul Technical University, Faculty of Science and Letters, Department of Mathematics,
34469 Maslak, Istanbul, Turkey}
\email{turgayn@itu.edu.tr}


\begin{abstract}
We define general rotational surfaces of elliptic and hyperbolic type in the pseudo-Euclidean 4-space with neutral metric  which are analogous to the general rotational surfaces of C. Moore in the Euclidean 4-space. We study Lorentz general rotational surfaces with plane  meridian curves and  give the complete classification of minimal general rotational surfaces of elliptic and hyperbolic type,   general rotational surfaces with parallel normalized mean curvature vector field, flat  general rotational surfaces, and general rotational surfaces with flat normal connection.
\end{abstract}

\maketitle

\section{Introduction}

 In  \cite{Moore1} C. Moore introduced general rotational surfaces in the 4-dimensional Euclidean space $\E^4$ and described a special case of general rotational surfaces with constant Gauss curvature \cite{Moore2}.
In \cite{Mil} the second author studied general rotational surfaces whose meridians lie in two-dimensional planes  and completely  classified
all minimal super-conformal general rotational surfaces in $\E^4$. The minimal non-super-conformal general rotational surfaces with plane meridian curves are described by U. Dursun and the third author in \cite{Dur&Tur2}.
In \cite{GM2}  the complete classification of general rotational surfaces consisting of parabolic points is given. General rotational surfaces with plane meridian curves and pointwise 1-type Gauss map are studied in \cite{Dur&Tur1}. 

In \cite{GM-TJM} G. Ganchev and the second author studied spacelike general rotational surfaces in the 4-dimensional Minkowski space $\E^4_1$  that are analogous to the general rotational surfaces of C. Moore in $\E^4$ and described analytically  flat general rotation surfaces and  general rotational surfaces with flat normal connection. The classification  of minimal general rotational surfaces in $\E^4_1$  and  general rotational surfaces consisting of parabolic points is also given in \cite{GM-TJM}. Spacelike general rotational surfaces in  $\E^4_1$  with meridian curves lying in 2-dimensional planes and having pointwise 1-type Gauss map  are studied in \cite{Dur}. 

In the present paper we define general rotational surfaces of elliptic and hyperbolic type in the pseudo-Euclidean 4-space with neutral metric $\E^4_2$ which are analogous to the general rotational surfaces in $\E^4$ and $\E^4_1$. We study Lorentz general rotational surfaces with meridian curves lying in 2-dimensional planes. In Theorem \ref{T:rotational-min-ell} and Theorem \ref{T:rotational-min-hyp} we give the complete classification of minimal general rotational surfaces of elliptic and hyperbolic type. Theorem \ref{T:Parallel-ell}  and Theorem \ref{T:Parallel-hyp} classify general rotational surfaces of elliptic and hyperbolic type with parallel normalized mean curvature vector field. The classification of flat  general rotational surfaces of elliptic and hyperbolic type is given in  Theorem \ref{T:rotational-flat-ell} and Theorem \ref{T:rotational-flat-hyp}, respectively. In the last section we describe all  general rotational surfaces of elliptic and hyperbolic type  with flat normal connection (Theorem \ref{T:flat normal} and Theorem \ref{T:flat normal-hyp}). 

\section{Preliminaries}

Let  $\E^4_2$ be the pseudo-Euclidean 4-space  endowed with the canonical pseudo-Euclidean metric of index 2 given by
$$g_0 = dx_1^2 + dx_2^2 - dx_3^2 - dx_4^2,$$
where $(x_1, x_2, x_3, x_4)$ is a rectangular coordinate system of $\E^4_2$. As usual, we denote by
$\langle .\, , . \rangle$ the indefinite inner scalar product with respect to $g_0$.
 A non-zero vector $v$ is called  \emph{spacelike} (respectively, \emph{timelike}) if $\langle v, v \rangle > 0$ (respectively, $\langle v, v \rangle < 0$).
 A vector $v$ is called \emph{lightlike} if it is nonzero and satisfies $\langle v, v \rangle = 0$.

A surface $M^2_1$ in $\E^4_2$ is called \emph{Lorentz}  if the
induced  metric $g$ on $M^2_1$ is Lorentzian, i.e. at each point $p
\in M^2_1$ we have the following decomposition
$$\E^4_2 = T_pM^2_1 \oplus N_pM^2_1$$
with the property that the restriction of the metric
onto the tangent space $T_pM^2_1$ is of
signature $(1,1)$, and the restriction of the metric onto the normal space $N_pM^2_1$ is of signature $(1,1)$.

Denote by $\nabla$ and $\nabla'$ the Levi Civita connections of $M^2_1$  and $\E^4_2$, respectively.
Let $x$ and $y$ be vector fields tangent to $M^2_1$ and $\xi$ be a normal vector field.
The formulas of Gauss and Weingarten are given respectively by
$$\begin{array}{l}
\vspace{2mm}
\nabla'_xy = \nabla_xy + \sigma(x,y);\\
\vspace{2mm}
\nabla'_x \xi = - A_{\xi} x + D_x \xi,
\end{array}$$
where $\sigma$ is the second fundamental form, $D$ is the normal
connection, and $A_{\xi}$ is the shape operator  with respect to
$\xi$. In general, $A_{\xi}$ is not diagonalizable.

The mean curvature vector  field $H$ of   $M^2_1$
is defined as $H = \frac{1}{2}\,  \tr\, \sigma$.
A  surface $M^2_1$  is called \emph{minimal} if its mean curvature vector vanishes identically, i.e. 
$H =0$.
A  surface $M^2_1$  is called  \emph{quasi-minimal} (or \textit{pseudo-minimal}) if its
mean curvature vector is lightlike at each point, i.e. $H \neq 0$ and $\langle H, H \rangle =0$  \cite{Rosca}.

A normal vector field $\xi$ on $M^2_1$ is called \emph{parallel in the normal bundle} (or simply \emph{parallel}) if $D{\xi}=0$ holds identically \cite{Chen}.
The surface $M^2_1$ is said to have \emph{parallel mean curvature vector field} if its mean curvature vector $H$ is parallel, i.e.
$D H =0$.

A natural extension of the class of surfaces with parallel mean curvature vector field are surfaces with parallel
normalized mean curvature vector field. A surface  $M^2_1$ is said to have \textit{parallel normalized mean curvature vector field} if 
the mean curvature vector $H$ is non-zero and  there exists a unit vector field in the direction of the mean curvature vector field
which is parallel in the normal bundle \cite{Chen-MM}.
It is easy to see  that if $M^2_1$ is a surface  with non-zero parallel mean curvature vector field $H$ (i.e. $DH = 0$),
then $M^2_1$ is a surface with parallel normalized mean curvature vector field, but the converse is not true in general.
It is true only in the case $\Vert H \Vert = const$.

\section{General Rotational Surfaces of Elliptic and Hyperbolic Type}

General rotational surfaces  in the Euclidean 4-space $\E^4$ were
introduced by  C. Moore \cite{Moore1} as follows. Let $m: x(u) = \left(
x^1(u), x^2(u),  x^3(u), x^4(u)\right)$; $ u \in J \subset \R$ be
a smooth curve in $\E^4$, and $\alpha$, $\beta$ be constants. A
general rotation of the meridian curve $m$  is defined by
$$X(u,v)= \left( X^1(u,v), X^2(u,v),  X^3(u,v), X^4(u,v)\right),$$
where
$$\begin{array}{ll}
\vspace{2mm} X^1(u,v) = x^1(u)\cos\alpha v - x^2(u)\sin\alpha v; &
\qquad
X^3(u,v) = x^3(u)\cos\beta v - x^4(u)\sin\beta v; \\
\vspace{2mm} X^2(u,v) = x^1(u)\sin\alpha v + x^2(u)\cos\alpha v;&
\qquad X^4(u,v) = x^3(u)\sin\beta v + x^4(u)\cos\beta v.
\end{array}$$
If $\beta = 0$, $x^2(u) = 0$  the plane $Oe_3e_4$ is
fixed and one gets the classical rotation about a fixed
two-dimensional axis.

Similarly to the general rotations in the Euclidean space $\E^4$ one can
consider general rotational surfaces  in the Minkowski 4-space
$\E^4_1$ (see \cite{GM-TJM} and \cite{Dur}).

Now we shall define general rotational surfaces of Moore type in the pseudo-Euclidean 4-space $\E^4_2$.

Let $Oe_1e_2e_3e_4$ be an orthonormal base of $\E^4_2$, such that $\langle e_1, e_1\rangle=\langle e_2, e_2\rangle = 1$, and
$\langle e_3, e_3\rangle=\langle e_4, e_4\rangle = -1$. 
Let $m: x(u) = \left( x^1(u), x^2(u),  x^3(u), x^4(u)\right)$; $ u \in J \subset \R$ be
a smooth spacelike or timelike curve in $\E^4_2$, and $\alpha$,
$\beta$ be constants. 
A general rotational surface of elliptic type can be defined as follows:
$$X(u,v)= \left( X^1(u,v), X^2(u,v),  X^3(u,v), X^4(u,v)\right),$$
where
\begin{equation} \label{E:Eq-gen-ell}
\begin{array}{ll}
\vspace{2mm} X^1(u,v) = x^1(u)\cos\alpha v - x^2(u)\sin\alpha v; &
\quad
X^3(u,v) = x^3(u)\cos \beta v - x^4(u)\sin \beta v; \\
\vspace{2mm} X^2(u,v) = x^1(u)\sin\alpha v + x^2(u)\cos\alpha v;&
\quad X^4(u,v) = x^3(u)\sin\beta v + x^4(u)\cos\beta v.
\end{array}
\end{equation}

If $\beta = 0$, $x^2(u) = 0$   one
gets the surface with parametrization 
$$X(u,v) = \left(x^1(u)\cos\alpha v, x^1(u)\sin\alpha v,  x^3(u), x^4(u)\right),$$  which is a standard rotation  of elliptic type about the 2-dimensional axis $Oe_3e_4$.
In the case $\alpha = 0$, $x^4(u) = 0$  we get the surface 
$$X(u,v) = \left(x^1(u),  x^2(u),  x^3(u)\cos \beta v, x^3(u)\sin \beta v\right),$$ which is a standard
 rotation of elliptic type about  $Oe_1e_2$.
 If $\alpha > 0$ and $\beta > 0$
the surface defined by \eqref{E:Eq-gen-ell} is analogous to the general rotational
surface of C. Moore in $\E^4$.

In the present paper we shall consider Lorentz general rotational surfaces of elliptic type for which $\alpha > 0$,  $\beta > 0$, $x^2(u) = x^4(u) = 0$. In this case the meridian curve $m$ lies in two-dimensional plane.  

Similarly to the general rotational surfaces of elliptic type we define 
 general rotational surfaces of hyperbolic type in $\E^4_2$ as follows:
\begin{equation} \label{E:Eq-gen-hyp}
\begin{array}{ll}
\vspace{2mm} X^1(u,v) = x^1(u)\cosh\alpha v + x^3(u)\sinh\alpha v; &
\quad
X^3(u,v) = x^1(u)\sinh \alpha v + x^3(u)\cosh \alpha v; \\
\vspace{2mm} X^2(u,v) = x^2(u)\cosh\beta  v + x^4(u)\sinh\beta v;&
\quad X^4(u,v) = x^2(u)\sinh\beta v + x^4(u)\cosh\beta v.
\end{array}
\end{equation}

In the case $\beta = 0$, $x^3(u) = 0$  we obtain 
the surface with parametrization 
$$X(u,v) = \left(x^1(u)\cosh\alpha v, x^2(u), x^1(u)\sinh\alpha v, x^4(u)\right),$$  which is a standard rotation  of hyperbolic type about $Oe_2e_4$.
In the case $\alpha = 0$, $x^4(u) = 0$  we get the surface 
$$X(u,v) = \left(x^1(u),  x^2(u)\cosh \beta v, x^3(u),  x^2(u)\sinh \beta v\right),$$ which is a standard
 rotation of hyperbolic type about  $Oe_1e_3$.
 If $\alpha > 0$ and $\beta > 0$ the surface defined by \eqref{E:Eq-gen-hyp} is a general rotational
surface of hyperbolic type in $\E^4_2$. We shall consider Lorentz general rotational surfaces of hyperbolic type for which  $\alpha > 0$,  $\beta > 0$, $x^3(u) = x^4(u) = 0$.

\vskip 2mm
\subsection{General rotational surfaces of elliptic type with plane meridians}

Let $\mathcal{M}_1$ be the general rotational surface of elliptic type  defined  by:
\begin{equation} \label{E:Eq-ell}
\mathcal{M}_1: z(u,v) = \left(f(u) \cos\alpha v, f(u) \sin \alpha v, g(u) \cos \beta v, g(u)\sin \beta v \right),
\end{equation}
where $u \in J \subset \R, \,\,  v \in [0;
2\pi)$, $f(u)$ and $g(u)$ are smooth functions satisfying
$\alpha^2 f^2(u) - \beta^2 g^2(u)
< 0 , \,\, f'\,^2(u) - g'\,^2(u) > 0$, and  $\alpha$, $\beta$ are positive constants. 

The tangent frame field $T_p \mathcal{M}_1$ is determined by the vector fields
\begin{equation*}
\begin{array}{l}
\vspace{2mm} z_u = \left(f'(u) \cos\alpha v, f'(u) \sin \alpha v, g'(u) \cos \beta v, g'(u)\sin \beta v \right);\\
\vspace{2mm} z_v = \left(-\alpha f(u) \sin\alpha v, \alpha f(u) \cos \alpha v, -\beta g(u) \sin \beta v, \beta g(u)\cos \beta v \right).\
\end{array}
\end{equation*}
The coefficients of the first fundamental form of  $\mathcal{M}_1$ are expressed by:
\begin{equation*}
\begin{array}{l}
\vspace{2mm} E =\langle z_u,z_u \rangle = f'^2(u) - g'^2(v) > 0;\\
\vspace{2mm} F = \langle z_u,z_v \rangle = 0;\\
\vspace{2mm} G = \langle z_v,z_v \rangle = \alpha^2 f^2(u) - \beta^2 g^2(u) < 0.
\end{array}
\end{equation*}
So, $\mathcal{M}_1$ is a Lorentz surface in $\E^4_2$.

We consider the following  orthonormal tangent frame field: 
\begin{equation*}
\vspace{2mm} x =\ds \frac{z_u}{\sqrt{E}}; \quad    y = \ds \frac{z_v}{\sqrt{-G}},
\end{equation*}
 which satisfy  $\langle x, x \rangle = 1$, $\langle y, y \rangle = -1$, and $\langle x, y \rangle = 0$. 
Let  $ n_1$ and $n_2$  be the normal vector fields defined by:
\begin{equation*}
\begin{array}{l}
\vspace{2mm} 
n_1 = \frac{1}{\sqrt{-\alpha^2 f^2(u) + \beta^2 g^2(u)}}\left(\beta g(u)\sin\alpha v, -\beta g(u)\cos\alpha v, \alpha f(u)\sin\beta v, -\alpha f(u)\cos\beta v\right);\\
\vspace{2mm} 
n_2 = \frac{1}{\sqrt{f'^2(u) - g'^2(v)}}\left(\ds{g'(u)\cos\alpha v}, \ds{g'(u)\sin\alpha v}, \ds{f'(u)\cos\beta v}, \ds{f'(u)\sin\beta v}\right).
\end{array}
\end{equation*}
Note that $\langle n_1,n_1 \rangle = 1$, $\langle n_2,n_2 \rangle = -1$, $\langle n_1,n_2 \rangle = 0$. 
Calculating the second derivatives 
\begin{equation*}
\begin{array}{l}
\vspace{2mm} z_{uu} = \left(f''(u) \cos\alpha v, f''(u) \sin \alpha v, g''(u) \cos \beta v, g''(u)\sin \beta v \right);\\
\vspace{2mm} z_{uv} = \left(-\alpha f'(u) \sin\alpha v, \alpha f'(u) \cos \alpha v, -\beta g'(u) \sin \beta v, \beta g'(u)\cos \beta v \right);\\
\vspace{2mm} z_{vv} = \left(-\alpha^2 f(u) \cos\alpha v, -\alpha^2 f(u) \sin \alpha v, -\beta^2 g(u) \cos \beta v, -\beta^2 g(u)\sin \beta v \right),
\end{array}
\end{equation*}
we obtain the following components of the second fundamental tensor:
$$\begin{array}{ll}
\vspace{2mm}
\langle z_{uu},n_{1} \rangle = 0; & \quad \langle z_{uu},n_{2} \rangle = \ds{\frac{f''(u) g'(u) - g''(u) f'(u)}{\sqrt{f'^2(u) - g'^2(u)}}}; \\
\vspace{2mm}
 \langle z_{uv},n_{1} \rangle = \ds{\frac{\alpha\beta (f(u) g'(u) - f'(u) g(u))}{\sqrt{\beta^2 g^2(u) - \alpha^2 f^2(u)}}}; & \quad \langle z_{uv},n_{2} \rangle = 0;\\
\vspace{2mm}
\langle z_{vv},n_{1} \rangle = 0; & \quad \langle z_{vv},n_{2} \rangle = \ds{\frac{\beta^2f'(u) g(u) - \alpha^2f(u) g'(u))}{\sqrt{f'^2(u) - g'^2(u)}}}.
\end{array}$$
The above formulas imply
\begin{equation} \label{E:Eq-1}
\begin{array}{l}
\vspace{2mm}
\sigma (x, x) = \ds{\frac{f'(u) g''(u) - g'(u) f''(u)}{(\sqrt{f'^2(u) - g'^2(u)})^3}}\;n_2;\\
\vspace{2mm}
\sigma (x, y) = \ds{\frac{\alpha\beta (f(u) g'(u) - f'(u) g(u))}{\sqrt{f'^2(u) - g'^2(u)} (\beta^2 g^2(u) - \alpha^2 f^2(u))}}\; n_1;\\
\vspace{2mm}
\sigma (y, y) = \ds{\frac{(\alpha^2 f(u) g'(u)-\beta^2 f'(u) g(u))}{\sqrt{f'^2(u) - g'^2(u)} (\beta^2 g^2(u) - \alpha^2 f^2(u))}}\; n_2.
\end{array}
\end{equation}

Using that the Gauss curvature $K$  is determined by the formula 

\begin{equation*}
\ds{K={\frac{\langle \sigma(x, x), \sigma(y, y) \rangle - \langle \sigma(x, y), \sigma(x, y) \rangle}{\langle x, x\rangle \langle y, y \rangle - \langle x, y \rangle^2}}},
\end{equation*}
we obtain 
\begin{equation*} \label{E:Eq-K}
\ds{K={\frac{\alpha^2\beta^2 (f'^2-g'^2) (f g' - f' g)^2 - (\beta^2 g^2 - \alpha^2 f^2) (\beta^2 f' g - \alpha^2 f g') (f' g'' - f'' g')}{(f'^2 - g'^2)^2 (\beta^2 g^2 - \alpha^2 f^2)^2}}}.
\end{equation*}

Calculating the derivatives of the normal vector fields $n_1$ and $n_2$ we get
\begin{equation} \label{E:Eq-2}
\begin{array}{ll}
\vspace{2mm}
D_x n_1 = 0; & \quad D_x n_2 = 0;\\
\vspace{2mm}
D_y n_1 =  \ds{\frac{\alpha \beta (f f' - g g')}{\sqrt{f'^2 - g'^2}(\beta^2 g^2 - \alpha^2 f^2)}} \,n_2; & \quad 
D_y n_2 =  \ds{\frac{\alpha \beta (f f' - g g')}{\sqrt{f'^2 - g'^2}(\beta^2 g^2 - \alpha^2 f^2)}}\, n_1.
\end{array}
\end{equation}
For the tangent vector fields $x$ and $y$ we obtain the following derivative formulas:
\begin{equation} \label{E:Eq-3}
\begin{array}{ll}
\vspace{2mm}
\nabla_x x = 0; & \quad \nabla_xy = 0;\\
\vspace{2mm}
\nabla_y x = - \frac{\alpha^2 f f' - \beta^2 g g'}{\sqrt{f'^2 - g'^2}(\beta^2 g^2 - \alpha^2 f^2)}\,y; & \quad 
\nabla_y y = - \frac{\alpha^2 f f' - \beta^2 g g'}{\sqrt{f'^2 - g'^2}(\beta^2 g^2 - \alpha^2 f^2)} \, x.
\end{array}
\end{equation}

The curvature of the normal connection $\varkappa$ of $\mathcal{M}_1$ is determined by 
\begin{equation*}
\varkappa = \frac{\langle R^{\bot}(x,y)n_1, n_2\rangle}{\langle x, x\rangle \langle y, y \rangle - \langle x, y \rangle^2} = - \langle D_xD_yn_1 - D_yD_xn_1 - D_{[x,y]} n_1, n_2 \rangle.
\end{equation*}
Formulas \eqref{E:Eq-2} and  \eqref{E:Eq-3} imply that the curvature of the normal connection is given by the following expression:
\begin{equation} \label{E:Eq-kappa}
\varkappa=\ds{\frac{-\alpha \beta (f g' - g f') \left((\beta^2  g^2 - \alpha^2 f^2)(g' f'' - f' g'') + (f'^2 -
g'^2) (\beta^2 g f' - \alpha^2 f g') \right)}{(f'^2 - g'^2)^2 (\beta^2 g^2 - \alpha^2 f^2)^2} }.
\end{equation}

The normal mean curvature vector field  $H$ is defined by the formula
\begin{equation*} 
\ds{H={\frac{\sigma(x, x) - \sigma(y, y)}{2}}}.
\end{equation*}
Using the expressions for  $\sigma(x, x)$ and $\sigma(y, y)$ in formulas \eqref{E:Eq-1} we get:
\begin{equation} \label{E:Eq-H}
\ds{H = \frac{(f'^2 - g'^2) (\beta^2 g f' - \alpha^2 f g') - (\beta^2 g^2 - \alpha^2 f^2)(f'' g' - f' g'')}{ 2 \left(f'^2 - g'^2 \right)^{\frac{3}{2}} (\beta^2 g^2 - \alpha^2 f^2)}\; n_2}.
\end{equation}

Formula \eqref{E:Eq-H} shows that in the case $H \neq 0$  the normalized  mean curvature vector field  of a general rotational surface of elliptic type with plane meridian curves is timelike. So, we can formulate the following statement:

\begin{prop} 
There are no quasi-minimal general rotational surfaces of elliptic type with plane meridian curves.
\end{prop}

\begin{rem} Let us note that in the class of standard rotational surfaces in $\E^4_2$ there exist quasi-minimal surfaces and they are described in  \cite{GM5}.
\end{rem}

In our further considerations we shall use the following notations:
\begin{equation} \label{E:Eq-inv} 
\begin{array}{ll}
\vspace{2mm}
\nu_1 = \ds{\frac{g' f'' - f' g''}{\left(f'^2 - g'^2\right)^{\frac{3}{2}}}}; & \qquad
\nu_2 = \ds{\frac{\beta^2 g f' - \alpha^2 f g'}{\sqrt{f'^2 - g'^2}(\beta^2 g^2 - \alpha^2 f^2)}};\\
\vspace{2mm}
\mu = \ds{\frac{\alpha \beta (f g' - g f')}{\sqrt{(f'^2 - g'^2)}(\beta^2 g^2 - \alpha^2 f^2)}}; & \qquad 
\gamma_2 = \ds{\frac{\alpha^2 f f' - \beta^2 g g'}{\sqrt{f'^2 - g'^2}(\beta^2 g^2 - \alpha^2 f^2)}};\\
\vspace{2mm}
\beta_2 = \ds{\frac{\alpha \beta (f f' - g g')}{\sqrt{f'^2 - g'^2}(\beta^2 g^2 - \alpha^2 f^2)}}.   & 
\end{array}
\end{equation}
Thus, equalities \eqref{E:Eq-1}, \eqref{E:Eq-2}, and \eqref{E:Eq-3} give  the following derivative formulas of  $\mathcal{M}_1$:
\begin{equation}\label{E:Eq-4}
\begin{array}{ll}
\vspace{2mm}\nabla'_{x}x = -\nu_{1}\,n_2; & \qquad \nabla'_{x}n_1 = \mu\,y;\\
\vspace{2mm}\nabla'_{x}y =  \mu\,n_1; & \qquad \nabla'_{y}n_1 = -\mu\,x + \beta_{2}\,n_2;\\
\vspace{2mm}\nabla'_{y}x = -\gamma_{2}\,y + \mu\,n_1; & \qquad \nabla'_{x} n_2 = -\nu_{1}\,x;\\
\vspace{2mm}\nabla'_{y}y = -\gamma_{2}\,x - \nu_{2}\,n_2; & \qquad \nabla'_{y} n_2 = \nu_{2}\,y + \beta_{2}\,n_1.
\end{array}
\end{equation}

\begin{rem} In  \cite{AGM}, an invariant local theory of Lorentz surfaces in the pseudo-Euclidean space $\E^4_2$ is developed  and a family  of eight  geometric functions is introduced. It is proved that these geometric functions determine the surface up to a rigid motion in  $\E^4_2$.  
The functions $\nu_1$, $\nu_2$, $\mu$,  $\gamma_2$,  $\beta_2$ given in \eqref{E:Eq-inv} are the geometric functions of the general rotational surface  $\mathcal{M}_1$ in the sense of \cite{AGM} (note that the other three geometric functions of $\mathcal{M}_1$ are equal to zero).
\end{rem}

\subsection{General rotational surfaces of hyperbolic type with plane meridians}

Now we shall consider   general rotational surfaces of hyperbolic  type with plane meridian curves. Let $\mathcal{M}_2$ be the surface defined by 
\begin{equation}  \label{E:Eq-hyp}
\mathcal{M}_2: z(u,v) = \left(f(u) \cosh\alpha v, g(u) \cosh \beta v, f(u) \sinh \alpha v, g(u)\sinh \beta v \right),
\end{equation}
where $u \in J \subset \R, \,\,  v \in [0;2\pi)$, $f(u)$ and $g(u)$ are smooth functions, which satisfy the conditions
$\alpha^2 f^2(u) + \beta^2 g^2(u) > 0 , \,\, f'^2(u) + g'^2(u) > 0$, $\alpha=const >0$, $\beta=const >0$.

The tangent vector fields $z_u$ and $z_v$ of  $\mathcal{M}_2$ are:
\begin{equation*}
\begin{array}{l}
\vspace{2mm} z_u =\left(f'(u) \cosh\alpha v, g'(u) \cosh \beta v, f'(u) \sinh \alpha v, g'(u)\sinh \beta v \right),\\
\vspace{2mm} z_v = \left(\alpha f(u) \sinh\alpha v, \beta g(u) \sinh \beta v, \alpha f(u) \cosh \alpha v, \beta g(u)\cosh \beta v \right).\
\end{array}
\end{equation*}
The coefficients of the first fundamental form of $\mathcal{M}_2$ are given by:
$$\vspace{2mm} E = f'^2(u) + g'^2(v),\ \quad
 F =  0, \quad  G =  -(\alpha^2 f^2 + \beta^2 g^2).$$
So,  $\mathcal{M}_2$ is a Lorentz surface.

 Now we consider the following  orthonormal tangent frame field 
$x =\ds \frac{z_u}{\sqrt{E}}, \;   y = \ds \frac{z_v}{\sqrt{-G}}$, which satisfies $\langle x,x \rangle = 1$, $\langle y,y \rangle = -1$, $\langle x,y \rangle = 0$. 
We choose the following normal vector  fields of $\mathcal{M}_2$:
\begin{equation*}
\begin{array}{l}
\vspace{2mm} n_1  = \frac{1} {\sqrt{f'^2(u) + g'^2(v)}} \left(\ds{g'(u)\cosh\alpha v}, \ds{-{f'(u)\cosh\beta v}, \ds{g'(u)\sinh\alpha v}}, \ds{-f'(u)\sinh\beta v}\right);\\
\vspace{2mm} n_2 = \frac{1}{\sqrt{\alpha^2 f^2(u) + \beta^2 g^2(u)}}\left(\ds{\beta g(u)\sinh\alpha v}, \ds{-\alpha f(u)\sinh\beta v}, \ds{\beta g(u)\cosh\alpha v}, \ds{-\alpha f(u)\cosh\beta v}\right),
\end{array}
\end{equation*}
for which we have  $\langle n_1,n_1 \rangle = 1$, $\langle n_2,n_2 \rangle = -1$, $\langle n_1,n_2 \rangle = 0$. 

With respect to the frame field $\{x,y,n_1,n_2\}$ introduced above we obtain the following derivative formulas of $\mathcal{M}_2$:
\begin{equation}\label{E:Eq-4-hyp}
\begin{array}{ll}
\vspace{2mm}\nabla'_{x}x = \nu_{1}\,n_1; & \qquad \nabla'_{x}n_1 = -\nu_{1}\,x;\\
\vspace{2mm}\nabla'_{x}y =  -\mu\,n_2; & \qquad \nabla'_{y}n_1 = \nu_{2}\,y  - \beta_{2}\,n_2;\\
\vspace{2mm}\nabla'_{y}x = -\gamma_{2}\,y - \mu\,n_2; & \qquad \nabla'_{x} n_2 = \mu\,y;\\
\vspace{2mm}\nabla'_{y}y = -\gamma_{2}\,x + \nu_{2}\,n_1; & \qquad \nabla'_{y} n_2 = -\mu\,x - \beta_{2}\,n_1,
\end{array}
\end{equation}
where the functions  $\nu_1$, $\nu_2$, $\mu$,  $\gamma_2$,  $\beta_2$ are expressed by:
\begin{equation}\label{E:Eq-inv-hyp} 
\begin{array}{ll}
\vspace{2mm} 
\nu_1 = \ds{\frac{f'' g' - f' g''}{\left(\sqrt{f'^2 + g'^2}\right)^3}}; & \qquad
\nu_2 = \ds{\frac{\alpha^2 f g' - \beta^2 g f'}{\sqrt{f'^2 + g'^2}(\alpha^2 f^2 + \beta^2 g^2)}};\\
\vspace{2mm}
\mu = \ds{\frac{\alpha \beta (fg' - f'g)}{\sqrt{f'^2 + g'^2}(\alpha^2 f^2 + \beta^2 g^2)}};
 & \qquad \gamma_2 = \ds{-\frac{\alpha^2 f f' + \beta^2 g g'}{\sqrt{f'^2 + g'^2}(\alpha^2 f^2 + \beta^2 g^2)}};\\
\vspace{2mm}
 \beta_2 = \ds{-\frac{\alpha \beta (f f' + g g')}{\sqrt{f'^2 + g'^2}(\alpha^2 f^2 + \beta^2 g^2)}}. &
\end{array}
\end{equation}
These five functions are the geometric functions of the general rotational surface of hyperbolic type  $\mathcal{M}_2$ in the sense of \cite{AGM}. The other three geometric functions of $\mathcal{M}_2$ are equal to zero.

\vskip 2mm
Similarly to the elliptic case we obtain the following expressions for the Gauss curvature $K$, the curvature of the normal connection $\varkappa$,  and the mean curvature vector field $H$ of  $\mathcal{M}_2$:
\begin{equation*} \label{E:Eq-K-hyp}
\ds{K={- \frac{\alpha^2\beta^2 (f g' - f' g)^2 (f'^2 + g'^2)+(\alpha^2 f g' - \beta^2 f' g) (f'' g' - f' g'') (\alpha^2 f^2 + \beta^2 g^2)}{(f'^2 + g'^2)^2 (\alpha^2 f^2 + \beta^2 g^2)^2}}};
\end{equation*}
\begin{equation*} \label{E:Eq-kappa-hyp}
\varkappa=\ds{\frac{\alpha \beta (f g' - f' g) \left((\alpha^2 f^2 + \beta^2  g^2)(f'' g' - f'  g'') + (f'^2 +
g'^2) (\alpha^2 f g' -  \beta^2 g f') \right)}{(f'^2 + g'^2)^2 (\alpha^2 f^2 + \beta^2 g^2)^2}};
\end{equation*}
\begin{equation*}\label{E:eq-H-hyp}
\ds{H = \frac{(f'^2 + g'^2) (\beta^2 f' g - \alpha^2 f g') + (\alpha^2 f^2 + \beta^2 g^2)(f'' g' - f' g'')}{ 2 \left(f'^2 + g'^2 \right)^{\frac{3}{2}} (\alpha^2 f^2 + \beta^2 g^2)}\; n_1}.
\end{equation*}

Note that in the case $H \neq 0$  the normalized  mean curvature vector field  of a general rotational surface of hyperbolic type with plane meridian curves is timelike. So, we can formulate the following statement:

\begin{prop} 
There are no quasi-minimal general rotational surfaces of hyperbolic type with plane meridian curves.
\end{prop}

\vskip 2mm
If  a Lorentz surface in $\E^4_2$ is neither minimal nor quasi-minimal, then it is worth considering the allied mean curvature vector field of the surface. The notion of allied mean curvature vector field of a Lorentz surface $M^2_1$ in $\E^4_2$ is defined by the formula
\begin{equation*} \label{E:eq17}
a(H)=\ds{\frac{\|H\|}{2} \, \tr(A_1 \circ A_2) \,n_2},
\end{equation*}
where
$\{n_{1}=\ds{\frac{H}{\|H\|}},n_{2}\}$ is an orthonormal base of the normal space of $M^2_1$;  
$A_1$ and $A_2$ are the shape operators corresponding to $n_1$ and $n_2$, respectively. $M^2_1$ is said to be a \textit{Chen surface} (\textit{Chen submanifold} in $\E^4_2$) if $a(H)$ vanishes identically \cite{Chen}. 

Using formulas \eqref{E:Eq-4}  and \eqref{E:Eq-4-hyp} for the general rotational surface $\mathcal{M}_1$ and $\mathcal{M}_2$, respectively,  we can easily calculate that $a(H)= 0$. So, we get the following 
 
\begin{prop} 
Each general rotational surface of elliptic or hyperbolic  type in $\E^4_2$  is a Chen surface.
\end{prop}

\vskip 2mm
General rotational surfaces with plane meridian curves and pointwise 1-type Gauss map in $\E^4_2$ are studied in \cite{AY} under the assumption  $\alpha =\beta =1$ and $K=0$. In \cite{BeCaDu}, the classification of general rotational surfaces  having zero mean curvature and pointwise 1-type Gauss map of second kind is given.

A special class of rotational surfaces with constant mean curvature in  $\E^4_2$ is studied in \cite{Liu-Liu}. These are general rotational surfaces  with plane meridian curves in the case  $\alpha =\beta =1$ and $f(u) = \varphi(u) \sinh u; \, g(u) = \varphi(u) \cosh u$ for some smooth function $\varphi(u)$.

\vskip 1mm
In what follows  we study some basic classes of general rotational surfaces of elliptic and hyperbolic type in the case of arbitrary constants $\alpha >0$ and $\beta >0$. We give the classification of minimal general rotational surfaces, general rotational surfaces with parallel normalized mean curvature vector field,  and flat general rotational surfaces. We describe analytically the class of general rotational surfaces with flat normal connection.

\section{Minimal general rotational surfaces of elliptic or hyperbolic type}

The study of minimal surfaces is one of the main topics in classical differential geometry. Recall that a  surface in $\E^4_2$ is minimal if the mean curvature vector field $H = 0$. Each plane is a trivial minimal surface. We study only surfaces which do not contain any open part of a plane.

In the next theorem we classify all minimal general rotational surfaces of elliptic type.

\begin{thm} \label{T:rotational-min-ell}
Let $\mathcal{M}_1$ be a general rotational surface of elliptic type, defined by \eqref{E:Eq-ell}. Then $\mathcal{M}_1$  is minimal if and only if the meridian curve $m$ is  determined by  one of the following:

\vskip 2mm
(i) $\ds{f= cg^{\pm \frac{\alpha}{\beta}}}, \; c = const, \, c  \neq 0, \, \alpha \neq \beta$;

\vskip 2mm
(ii) $\ds{\arcsin\left(\frac{\alpha f}{\sqrt{A}}\right) = \pm \frac{\alpha}{\beta}\arcsin\left(\frac{\beta g}{\sqrt{A}}\right)} + C$, where $C = const, \, A = const,\,A > 0,  \, \alpha \neq \beta$;

\vskip 2mm
(iii) $\left(f + g\right)^2 = a\left(f - g\right)^2+ b$,
where $ a = const, \, a  \neq 0, \, b = const$, and $\alpha = \beta$.
\end{thm}

\noindent \emph{Proof:}  The mean curvature vector field $H$ of a general rotational surface of elliptic type is given by \eqref{E:Eq-H}. 
Hence, the surface $\mathcal{M}_1$ is minimal if and only if 
\begin{equation*}
\ds{\frac{g' f'' - f' g''}{f'^2 - g'^2} = \frac{\beta^2 g f' - \alpha^2 f g'}{\beta^2g^2 - \alpha^2f^2}}.
\end{equation*}
On the other hand, from \eqref{E:Eq-4} it follows that the mean curvature vector field $H$ is expressed
as $H = \ds{\frac{\nu_2 - \nu_1}{2}\; n_2}$,
where $\nu_1$ and $\nu_2$ are given in \eqref{E:Eq-inv}.
So, the condition $H=0$ is equivalent to $\nu_1=\nu_2$. 

Now, let  $\mathcal{M}_1$  is minimal, i.e. the equality $\nu_1=\nu_2$ holds.
Taking into account that $R'=0$ and using derivative formulas \eqref{E:Eq-4} we obtain  the equalities
\begin{equation} \label{E:eq-integ}
\begin{array}{l}
\vspace{2mm}
2 \mu \gamma_2 - \nu_1 \beta_2 = x(\mu);\\
\vspace{2mm}
\mu \beta_2 - 2 \nu_1 \gamma_2 = -x(\nu_1).
\end{array}
\end{equation}

 If we assume that $\mu = 0$,   from \eqref{E:Eq-inv}  we get  $f g' -g f' =0$, which implies $f = c g$, $c = const$.
Straightforward computations show that in this case  $\nu_1 = \nu_2 = 0$ and hence $\sigma(x, x) = \sigma(x, y) = \sigma(y, y)  = 0$. So,  the surface is totally geodesic, i.e.   $\mathcal{M}_1$ is locally a plane. 

So, further we consider $\mu \neq 0$. If we assume that $\nu_1 = \nu_2 = 0$, from  \eqref{E:Eq-inv} we get  $g' f'' - f' g'' =0$, which implies
\begin{equation*}
\ds{\frac{f''}{f'} = \frac{g''}{g'}}.
\end{equation*}
The solution of this equation is given by $f = a g +b$, $a = const \neq 0$, $b = const$. Without loss of generality we can consider $f(u)$ and $g(u)$ as follows: 
\begin{equation} \label{E:eq - f}
f(u) = a u + b; \quad  g(u)= u.
\end{equation}
In this case the function $\mu$ is expressed as
$\mu = \ds{\frac{\alpha\beta\,b}{\sqrt{(a^2-1)}(\beta^2u^2 - \alpha^2(au+b)^2)}}$.
On the other hand, using  \eqref{E:eq-integ} in the case  $\nu_1 =0$, we obtain $\mu \,\beta_2 =0$.
Since  $\mu \neq 0$, we get $\beta_2 = 0$. The  expression of $\beta_2$  in \eqref{E:Eq-inv} implies
$ff' - gg' = 0$. The last equality together with 
\eqref{E:eq - f} gives $a^2-1=0$,  which contradicts the assumption $f'^2 -g'^2 >0$.

So, further we consider  $\mu \neq 0$ and  $\nu_1 = \nu_2 \neq 0$.

In the case  $\alpha = \beta$, taking into account that $\nu_1 = \nu_2$, from \eqref{E:Eq-inv} we get the equation:
\begin{equation} \label{E:eq-f}
\ds{\frac{f''g'- f'g''}{f'^2 - g'^2}} = \ds{\frac{fg'- f'g}{f^2 - g^2}}.
\end{equation}
We denote  $\phi(u) = \ds{\frac{f(u)}{g(u)}}$ and $\psi(u) = \ds{\frac{f'(u)}{g'(u)}}$. Then equation  \eqref{E:eq-f} is written in the form
\begin{equation*}
\ds{\frac{\phi'(u)}{1 - \phi^2(u)}} =  \ds{\frac{\psi'(u)}{\psi^2(u) - 1}},
\end{equation*}
which implies the following equality
\begin{equation*}
\ds{\ln\left|\frac{1 + \phi}{1 - \phi}\right|} =  \ds{\ln\left|\frac{\psi - 1}{\psi + 1}\right|} + const.
\end{equation*}
Hence, we obtain
\begin{equation*}
\ds{\frac{(1 + \phi)(1 + \psi)}{(1 - \phi)(1 - \psi)}} =  const.
\end{equation*}
Now, using the expressions of $\phi(u)$ and $\psi(u)$ we get
\begin{equation*}
\left((f + g)^2\right)' = a\left((f - g)^2\right)', \quad  a = const \neq 0.
\end{equation*}
Hence, in the case  $\alpha = \beta$ the relation  between the functions  $f$ and $g$ is:
\begin{equation*}
\left(f + g\right)^2 = a\left(f - g\right)^2+ b,
\end{equation*}
where $ a = const \neq 0$,  $ b = const$. Thus we obtain case (iii) in the statement of the theorem.

Further we consider the case $\alpha \neq \beta$. 
If  $\mu^2 - \nu^2 = 0$, or equivalently $\nu_1 = \nu_2 = \pm \mu$, then the expressions of $\nu_2$ and $\mu$ in \eqref{E:Eq-inv} give the equation:
\begin{equation*}
\beta^2 g f' - \alpha^2 f g' = \pm \alpha \beta (f g' - g f'),
\end{equation*}
which can be written in the form:
\begin{equation*}
\frac{f'}{f} = \pm \frac{\alpha}{\beta}\, \frac{g'}{g} .
\end{equation*}
The solution of the last equation is:
\begin{equation*}
\ds{f= c g^{\pm \frac{\alpha}{\beta}}}, \quad c = const,\;  c \neq 0,
\end{equation*}
which gives    (i) in the statement of the theorem.

Now we consider the case $\mu^2 - \nu^2 \neq 0$. Denote $\nu:= \nu_1 = \nu_2$. Using equalities \eqref{E:eq-integ} we get the equations:
\begin{equation*}
x(\mu) = 2\mu\gamma_2 - \nu\beta_2;
\end{equation*}
\begin{equation*}
x(\nu) = 2\nu\gamma_2 - \mu\beta_2,
\end{equation*}
which imply 
\begin{equation}\label{E:eq-mu}
x(\mu^2 -\nu^2) = 4(\mu^2 - \nu^2)\gamma_2.
\end{equation}
Formula  \eqref{E:eq-mu} together with the equality  $\gamma_2 = -x(\ln(\sqrt{-G})$ give us  
\begin{equation*}
x(\ln |\mu^2 -\nu^2|) = -4x(\ln(\sqrt{-G}),
\end{equation*}
which implies 
$$x(\ln |(\mu^2 -\nu^2) G^2|) = 0.$$
Since the functions $\mu$, $\nu$, and $G$  do not depend on the parameter $v$, we obtain
\begin{equation*}
|(\mu^2 -\nu^2)G^2| =c^2, \quad   c = const \neq 0.
\end{equation*}
Using the expressions of $\mu$ and $\nu$ in \eqref{E:Eq-inv} and taking into account $f'^2 - g'^2 = 1$, we obtain the equation:
\begin{equation*}
\alpha^2\beta^2(fg'- f'g)^2 - (\beta^2f'g - \alpha^2fg')^2= c^2,
\end{equation*}
or equivalently:
\begin{equation*}
(\alpha^2 - \beta^2)(\beta^2f'^2g^2 - \alpha^2f^2(f'^2 - 1)) = c^2.
\end{equation*}
Since  $\alpha^2 - \beta^2 \neq 0$, putting $A=\frac{c^2}{\alpha^2-\beta^2}$, we obtain the following  expression for $f '^2$:
\begin{equation} \label{E:eq-f1}
f'^2 = \ds{\frac{A - \alpha^2f^2}{\beta^2g^2 - \alpha^2f^2}}.
\end{equation}
Using that $ g'^2 =f'^2 - 1$ we get:
\begin{equation} \label{E:eq-g1}
g'^2 = \ds{\frac{A - \beta^2g^2}{\beta^2g^2 - \alpha^2f^2}}.
\end{equation}
Since $f'^2 >0$ and $g'^2>0$, we get $A >\alpha^2 f^2$ and $A >\beta^2 g^2$, so $A>0$, i.e. $\alpha^2>\beta^2$.
Equations \eqref{E:eq-f1} and \eqref{E:eq-g1} imply
\begin{equation*}
\ds{\frac{f'^2}{A - \alpha^2f^2}} = \ds{\frac{g'^2}{A - \beta^2g^2}}.
\end{equation*}
The last equality gives the following relation between the functions $f$ and $g$:
\begin{equation*}
\ds{\arcsin\left(\frac{\alpha f}{\sqrt{A}}\right)} = \ds{\pm \frac{\alpha}{\beta}\arcsin\left(\frac{\beta g}{\sqrt{A}}\right)} + C, \quad C = const.
\end{equation*}
This corresponds to case (ii) in the statement of the theorem.

Conversely, if one of  (i),  (ii), or (iii) holds, then by straightforward computations it follows that $\nu_1=\nu_2$, i.e. $\mathcal{M}_1$ is a  minimal surface.
\qed

\vskip 3mm
The following theorem gives the classification of all minimal general rotational surfaces of hyperbolic type.

\begin{thm} \label{T:rotational-min-hyp} Let $\mathcal{M}_2$ be a general rotational surface of hyperbolic type, defined by \eqref{E:Eq-hyp}. Then $\mathcal{M}_2$ is minimal if and only if   the meridian curve $m$ is determined by one of the following:

\vskip 2mm
(i) $f = cg^{\mp \frac{\alpha}{\beta}}, \quad c = const, \, c \neq 0, \, \alpha \neq \beta$;

\vskip 2mm
(ii) $\alpha f + \sqrt{\alpha^2f^2 - A} = C\left(\beta g + \sqrt{\beta^2g^2 + A}\right)^{\ds{ \pm \frac{\alpha}{\beta}}}
$,  $C = const, A = const, A C \neq 0, \alpha \neq \beta$;

\vskip 2mm
(iii) $\arctan{\left(\ds{\frac{f'}{g'}}\right)} = -\arctan{\left(\ds{\frac{f}{g}}\right)} + c, \quad c = const$, and $\alpha = \beta$.
\end{thm}

\noindent \emph{Proof:}
Similarly to the proof of  Theorem \ref{T:rotational-min-ell}, we obtain that $\mu \neq 0$ and  $\nu_1 = \nu_2 \neq 0$. 
We denote $\nu:=\nu_1 = \nu_2$. 

In the case $\alpha = \beta$, using that  $\nu_1 = \nu_2$, from \eqref{E:Eq-inv-hyp} we get the equation 
\begin{equation*}\label{E:eq-f3}
\ds{\frac{f''g'- f'g''}{f'^2 + g'^2}} = \ds{\frac{fg'- f'g}{f^2 + g^2}},
\end{equation*} 
which implies 
\begin{equation*}
\left(\arctan{\left(\frac{f'}{g'}\right)}\right)' = - \left(\arctan{\left(\frac{f}{g}\right)}\right)'.
\end{equation*}
So,  the relation  between the functions $f$ and $g$ is given by:
\begin{equation*}
\arctan{\left(\ds{\frac{f'}{g'}}\right)} = -\arctan{\left(\ds{\frac{f}{g}}\right)} + c, \quad c = const,
\end{equation*}
which corresponds to case (iii) of the theorem.

Now, we consider the case  $\alpha \neq \beta$. If $\mu^2 - \nu^2 = 0$, then as in the elliptic case we get the equation:
\begin{equation*}
\alpha^2 f g' - \beta^2 g f' = \pm \alpha \beta (f g'- g f'),
\end{equation*}
whose solution is:
\begin{equation*}
\ds{f= c g^{\mp \frac{\alpha}{\beta}}}, \quad c = const, c \neq 0.
\end{equation*}
This is case (i) in the statement of the theorem.

If $\mu^2 - \nu^2 \neq 0$, then from $R'=0$ and \eqref{E:Eq-4-hyp} we obtain the equation 
\begin{equation*}\label{E:eq-mu-a}
x(\mu^2 -\nu^2) = 4(\mu^2 - \nu^2)\gamma_2.
\end{equation*}
On the other hand, we have $\gamma_2 = -x(\ln(\sqrt{-G})$, so  we get
\begin{equation*}
|\mu^2 -\nu^2|G^2 = c^2, \quad  c = const \neq 0.
\end{equation*}
Without loss of generality we assume that  $f'^2 + g'^2 = 1$. So, we obtain  the equation 
\begin{equation*}
(\beta^2f'^2g^2 - \alpha^2f^2(1 - f'^2))(\alpha^2- \beta^2) = c^2.
\end{equation*}
Since $\alpha^2 - \beta^2 \neq 0$, putting  $A =\ds{\frac{c^2}{\alpha^2- \beta^2}}$, we get the following expressions for $f'^2$ and $g'^2$:
\begin{equation*}
f'^2 = \ds{\frac{\alpha^2f^2 - A}{\alpha^2 f^2 + \beta^2 g^2}}, \quad g'^2 = \ds{\frac{\beta^2g^2 + A}{\alpha^2 f^2 + \beta^2 g^2}}, 
\end{equation*}
which imply
\begin{equation} \label{E:Eq-10}
\ds{\frac{f'^2}{\alpha^2f^2 - A}} = \ds{\frac{g'^2}{\beta^2g^2 + A}}.
\end{equation}
Note that $\alpha^2f^2 - A>0$ and $\beta^2g^2 + A>0$, since $f'^2>0$ and $g'^2>0$.
Equation \eqref{E:Eq-10} implies that the functions $f$ and $g$, which describe the meridian curve $m$ of $\mathcal{M}_2$, satisfy the following relation:
\begin{equation*}
\ds{\ln\left|\alpha f + \sqrt{\alpha^2f^2 - A}\right|} = \ds{\pm \frac{\alpha}{\beta}\ln\left|\beta g + \sqrt{\beta^2g^2 + A}\right|} + c,
\end{equation*}
 where $c = const$.
Finally, we obtain
\begin{equation*}
\alpha f + \sqrt{\alpha^2f^2 - A} = C\left(\beta g + \sqrt{\beta^2g^2 + A}\right)^{\ds{\pm \frac{\alpha}{\beta}}}.
\end{equation*}
This is case (ii) in the theorem.

\qed

\section{General rotational surfaces of elliptic or hyperbolic  type with parallel normalized mean curvature vector field} \label{S:parallel}

In the case $H \neq 0$ the normalized mean curvature vector field of the rotational  surface of elliptic type $\mathcal{M}_1$  is $n_2$. 
So, $\mathcal{M}_1$ has parallel normalized mean curvature vector field if and only if $D_x n_2 = D_y n_2 = 0$.
The next theorem describes all general rotational surfaces of elliptic type with parallel normalized mean curvature vector field.

\begin{thm} \label{T:Parallel-ell} Let $\mathcal{M}_1$  be a general rotational surface of elliptic type, defined by \eqref{E:Eq-ell}.  Then $\mathcal{M}_1$  has parallel normalized mean curvature vector field  if and only if the meridian curve $m$ is determined by 
\begin{equation*}
f(u) = \pm \sqrt{u^2 - C^2};\quad
g(u) = u, \quad C = const \neq 0.
\end{equation*}

\end{thm}

\noindent \emph{Proof:}
It follows from \eqref{E:Eq-2} that $D_x n_2 = D_y n_2 = 0$ if and only if the functions $f$ and $g$ satisfy the following differential equation:
$$f f'- g g'= 0,$$ 
which implies that $f^2 = g^2 + C_1$ for some constant $C_1$. Without loss of generality we can assume that $g(u) = u$. Then  $f(u) = \pm \sqrt{u^2 + C_1}$. Since $f'^2 - g'^2 >0$, we obtain $C_1 <0$. Hence $f(u) = \pm \sqrt{u^2 - C^2}$ for some constant $C \neq 0$.

\qed

\vskip 3mm

Using the last theorem and equality \eqref{E:Eq-H} we obtain that the mean curvature vector field $H$ of a general rotational surface of elliptic type with parallel normalized mean curvature vector field is given by:
\begin{equation*}
H = \ds\frac{1}{C} \, n_2.
\end{equation*}
Hence,  the condition on $\mathcal{M}_1$ to have  parallel normalized mean curvature vector field implies  $\langle H, H \rangle = - \frac{1}{C^2} = const$. So, the following statement holds: 

\begin {cor} \label{C:cor-parallel-H}
 If $\mathcal{M}_1$  has parallel normalized mean curvature vector field, then it has parallel mean curvature vector field. 
\end{cor}

\vskip 3mm
 The classification of  general rotational surfaces of hyperbolic type with parallel normalized mean curvature vector field is given in the next theorem.

\begin{thm} \label{T:Parallel-hyp}  Let $\mathcal{M}_2$  be a general rotational surface of hyperbolic type, defined by \eqref{E:Eq-hyp}. Then $\mathcal{M}_2$  has parallel normalized mean curvature vector field  if and only if the meridian curve $m$ is determined by 
\begin{equation*}
f(u) = \pm \sqrt{C^2 - u^2};\quad  
g(u) = u,  \quad C = const \neq 0.
\end{equation*}
\end{thm}

\noindent \emph{Proof:}
In the case $H \neq 0$ the normalized mean curvature vector field of  $\mathcal{M}_2$   is $n_1$. Using \eqref{E:Eq-4-hyp} and  \eqref{E:Eq-inv-hyp}, we get that $\mathcal{M}_2$ has parallel normalized mean curvature vector field if and only if the functions $f$ and $g$ satisfy the equation
$$ff'+ gg'= 0,$$ 
which implies $f^2  + g^2 = C^2$ for some constant $C \neq 0$. Without loss of generality we can assume that $g(u) = u$ and then the function $f(u)$ has the form
 $f(u) = \pm \sqrt{C^2 - u^2}$.

\qed

\vskip 2mm
Calculating the mean curvature vector field $H$ of a general rotational surface of hyperbolic type with parallel normalized mean curvature vector field, we obtain
\begin{equation*}
H = \ds{\mp \frac{1}{C}\, n_1}.
\end{equation*}

The last equality shows that $\| H \| = const$, so the following statement holds true:

\begin{cor} \label{C:cor-parallel-H-hyp}
If $\mathcal{M}_2$  has parallel normalized mean curvature vector field, then it has parallel mean curvature vector field. 
\end{cor}

\section{Flat general rotational surfaces of elliptic or hyperbolic type}
\vskip 3mm

In this section we give the classification of flat general rotational surfaces of elliptic or hyperbolic type.

\begin{thm} \label{T:rotational-flat-ell}
Let $\mathcal{M}_1$ be a general rotational surface of elliptic type, defined by \eqref{E:Eq-ell}. Then $\mathcal{M}_1$ is flat if and only if  the meridian curve $m$ is determined by one of the following:

\vskip 2mm
(i) $\beta^2g^2 - \alpha^2f^2 = a^2 (u + c)^2$,
where $a = const \neq 0$, $c = const$;

\vskip 2mm
(ii) $\alpha^2f^2 -\beta^2g^2 = C$, where $C = const, \, C <0$.
\end{thm}

\noindent \emph{Proof:} 
A Lorentz surface  is flat, if its Gauss curvature is zero. In terms of the functions $\mu$, $\nu_1$ and $\nu_2$ participating in formulas \eqref{E:Eq-4}, the condition on  $\mathcal{M}_1$  to have zero Gauss curvature is expressed as $\mu^2 + \nu_1\nu_2 =0$.

Now, let  $\mu^2 + \nu_1\nu_2 =0$, i.e. $K=0$. Using that $R'=0$ from derivative formulas \eqref{E:Eq-4} we obtain  
\begin{equation}{\label{E:eq-0}}
x(\gamma_2) = (\gamma_2)^2.
\end{equation}
Without loss of generality we can assume that $f'^2 - g'^2 = 1$. Then $x = \frac{\partial }{\partial u}$.

If $\gamma_2 \neq 0$ from equality  \eqref{E:eq-0} we get the equation :
\begin{equation*}
\ds{\frac{\gamma_2'}{(\gamma_2)^2} = 1},
\end{equation*}
whose  solution is
\begin{equation*}
\ds{\gamma_2(u)= - \frac{1}{u + c}}, \quad c =const.
\end{equation*}
On the other hand, $\gamma_2 = - x(\ln \sqrt{\beta^2g^2 - \alpha^2f^2})$. Consequently,
\begin{equation*}
(\ln \sqrt{\beta^2g^2 - \alpha^2f^2})'= \frac{1}{u + c}.
\end{equation*}
The last equation gives the following relation between the functions $f$ and $g$ determining the meridian curve $m$:
\begin{equation*}
\beta^2g^2 - \alpha^2f^2 = a^2 (u + c)^2,
\end{equation*}
where $a = const \neq 0$ and $c = const$.

If $\gamma_2 = 0$, from \eqref{E:Eq-inv} we get the equation
\begin{equation*}
\alpha^2 ff' -\beta^2 gg'=0,
\end{equation*}
which implies
\begin{equation*}
\alpha^2f^2 -\beta^2g^2 = C, \quad C = const.
\end{equation*}
Since we consider surfaces with $\alpha^2f^2 -\beta^2g^2 <0$, so the constant $C$ is negative.

Conversely, if  (i) or (ii) holds, then by straightforward computations we get  $\mu^2 + \nu_1\nu_2 =0$, i.e. the surface $\mathcal{M}_1$ is flat.

\qed

\vskip 3mm
The next theorem describes flat  general rotational surfaces of hyperbolic type.

\begin{thm} \label{T:rotational-flat-hyp} 
Let $\mathcal{M}_2$ be a general rotational surface of hyperbolic type, defined by \eqref{E:Eq-hyp}. Then $\mathcal{M}_2$ is flat if and only if the meridian curve $m$ is determined by one of the following:

\vskip 2mm
(i) $ \alpha^2f^2 + \beta^2g^2 = a^2 (u + c)^2, \quad 
a = const \neq 0$, $c = const$;

\vskip 2mm
(ii) $\alpha^2f^2 + \beta^2g^2 = C, \quad C = const$.
\end{thm}

\noindent \emph{Proof:} 
The proof in the hyperbolic case is similar to the elliptic one. Again we have that a general rotational surface of hyperbolic type is flat if and only if the functions $\mu$, $\nu_1$ and $\nu_2$ participating in \eqref{E:Eq-4-hyp} satisfy the condition $\mu^2 + \nu_1\nu_2 = 0$. The function $\gamma_2$ satisfies the equation 
$$x(\gamma_2) = (\gamma_2)^2.$$
In the case $\gamma_2 \neq 0$, assuming that  $f'^2 + g'^2 =1$,  we 
 obtain:
\begin{equation*}
\ds{\gamma_2(u)= - \frac{1}{u + c}}, \quad c =  const.
\end{equation*}
On the other hand, $\gamma_2 = - x(\ln \sqrt{\alpha^2f^2 + \beta^2g^2})$, so 
\begin{equation*}
\ds{(\ln \sqrt{\alpha^2f^2 + \beta^2g^2})'= \frac{1}{u + c}}.
\end{equation*}
The solution of the last differential equation is:

\begin{equation*}
\alpha^2f^2 + \beta^2g^2 = a^2(u + c)^2,
\end{equation*}
where $a = const \neq 0$, $c = const$.

If $\gamma_2 = 0$ we obtain
\begin{equation*}
\alpha^2 ff' + \beta^2 gg'=0,
\end{equation*}
or equivalently
\begin{equation*}
\alpha^2f^2 + \beta^2g^2 = C, \quad C = const.
\end{equation*}

\qed

\section{General rotational surfaces of elliptic or hyperbolic type with flat normal connection}
\vskip 3mm

A surface is said to have  flat normal connection if the curvature of the normal connection is zero.
Each surface with parallel normalized mean curvature vector field has flat normal connection. In Section \ref{S:parallel}, we classified all  general rotational surfaces with parallel normalized mean curvature vector field. Note that according to Corollary \ref{C:cor-parallel-H} and Corollary \ref{C:cor-parallel-H-hyp}, in the class of general rotational surfaces of both elliptic and hyperbolic type the condition of parallel  normalized mean curvature vector field is equivalent to the condition of parallel  mean curvature vector field. 
So, here we shall consider general rotational surfaces with flat normal connection but non-parallel mean curvature vector field, i.e. we assume that $\beta_2 \neq 0$.

\begin{thm}\label{T:flat normal}
Let $\mathcal{M}_1$ be a general rotational surface of elliptic type, defined by \eqref{E:Eq-ell}. Then $\mathcal{M}_1$ has flat normal connection and  non-parallel mean curvature vector field if and only if the meridian curve $m$ is determined by one of the following:

(i)  $f = c \,g$, where  $c = const, \, 1< c^2 < \frac{\beta^2}{\alpha^2}$, and  $\alpha< \beta$. In this case, $\mathcal{M}_1$  is a developable ruled surface in $\E^4_2$;

\vskip 1mm
(ii) $\frac{f f' - g g'}{ \sqrt{f'^2 - g'^2}\sqrt{\beta^2  g^2 - \alpha^2 f^2}} = C, \;  where \; C = const, \, C \neq 0$.
\end{thm}

\noindent \emph{Proof:} 
The curvature of the normal connection of the general rotational surface  $\mathcal{M}_1$ is
given by formula  \eqref{E:Eq-kappa}. Hence, $\mathcal{M}_1$ has flat normal connection if and only if
\begin{equation*} 
(f g' - g f') \left((\beta^2  g^2 - \alpha^2 f^2)(g' f'' - f' g'') + (f'^2 -
g'^2) (\beta^2 g f' - \alpha^2 f g') \right) = 0.
\end{equation*}
The last equation leads to the following two cases:

(i) $f g' - g f' = 0$, i.e. $f = c \,g$ for some constant $c \neq 0$. Since we consider surfaces for which $f'^2 - g'^2>0$ and 
$\alpha^2 f^2 - \beta^2  g^2 <0$, so the constant $c$ satisfies $1< c^2 < \frac{\beta^2}{\alpha^2}$, and hence $\alpha< \beta$.
By straightforward computations we get $\nu_1 =0$, $\nu_2 = \frac{c (\beta^2 - \alpha^2)}{g \sqrt{c^2-1}(\beta^2 - c^2 \alpha^2)} \neq 0$, $\mu =0$, 
$\gamma_2 = - \frac{1}{g\sqrt{c^2-1}} \neq 0$, $\beta_2 = \frac{\alpha \beta \sqrt{c^2-1}}{g (\beta^2 - c^2 \alpha^2)} \neq 0$. 
In this case  $\mathcal{M}_1$  is a ruled surface, since the meridian curve $m$ is a straight line.
Using formulas \eqref{E:Eq-4} we obtain that $\nabla'_x n_1 = 0$, $\nabla'_x n_2 = 0$, so,  the normal space is constant at each point of a  fixed straight line. Consequently,  $\mathcal{M}_1$ is a developable ruled surface in $\E^4_2$.

(ii) $f g' - g f' \neq 0$. In this case the functions $f$ and $g$ satisfy the  equation 
\begin{equation*} 
(\beta^2  g^2 - \alpha^2 f^2)(g' f'' - f' g'') + (f'^2 -
g'^2) (\beta^2 g f' - \alpha^2 f g')  = 0,
\end{equation*}
which is equivalent to $\nu_1 +\nu_2 = 0$.
 Using that $R'=0$, from derivative formulas \eqref{E:Eq-4} under the assumption $\nu_1 +\nu_2 = 0$, we obtain  $\gamma_2 \beta_2 = x(\beta_2)$.
Since $\beta_2 \neq 0$, we get the equation 
\begin{equation*}
\frac{x(\beta_2)}{\beta_2} = \gamma_2,
\end{equation*}
which implies $x(\ln |\beta_2|)= \gamma_2$. 
On the other hand, $\gamma_2 = - x(\ln \sqrt{-G})$. Consequently,
\begin{equation} \label{E:Eq-11}
x \left(\ln | \beta_2 \sqrt{-G} |\right)= 0.
\end{equation}
Now, using that the functions $\beta_2$  and $G$  do not depend on the parameter $v$, from \eqref{E:Eq-11} it follows that
\begin{equation*} 
\beta_2 \sqrt{-G} = const.
\end{equation*}
Having in mind the expression of $\beta_2$ given in \eqref{E:Eq-inv},  
we obtain that the functions $f$ and $g$ satisfy the following differential equation
\begin{equation*} \label{E:Eq-12}
\frac{f f' - g g'}{ \sqrt{f'^2 - g'^2}\sqrt{\beta^2  g^2 - \alpha^2 f^2}} = C,
\end{equation*}
where $C = const, \, C \neq 0$.

\qed

\vskip 3mm
Similarly to the proof of Theorem \ref{T:flat normal} we obtain the following characterization of general rotational surfaces of hyperbolic type with flat normal connection.

\begin{thm}\label{T:flat normal-hyp}
Let $\mathcal{M}_2$ be a general rotational surface of hyperbolic type, defined by \eqref{E:Eq-hyp}. Then $\mathcal{M}_2$ has flat normal connection  and  non-parallel mean curvature vector field if and only if the meridian curve $m$ is determined by one of the following:

(i)  $f = c \,g$, where  $c = const, \, c \neq 0$, and  $\alpha \neq \beta$. In this case, $\mathcal{M}_1$  is a developable ruled surface in $\E^4_2$;

\vskip 1mm
(ii) $\frac{f f' + g g'}{ \sqrt{f'^2 + g'^2}\sqrt{\alpha^2 f^2 + \beta^2  g^2}} = C, \;  where \; C = const, \, C \neq 0$.
\end{thm}

\vskip 5mm

\vskip 5mm \textbf{Acknowledgments:}
The second author is partially supported by the Bulgarian National Science Fund,
Ministry of Education and Science of Bulgaria under contract
DFNI-I 02/14. The third author is supported by T\"UB\.ITAK  (Project Name: Y\_EUCL2TIP, Project Number: 114F199).

This work was done during the third author's  visit at the Institute of Mathematics and Informatics, Bulgarian Academy of Sciences in June 2015.


\begin{thebibliography}{99}

\bibitem{AY}
Aksoyak F., Yayli Y., \textit{General rotational   surfaces with pointwise 1-type Gauss map in pseudo-Euclidean space $\E^4_2$}. Indian J. Pure Appl.  Math., \textbf{46}, no. 1  (2015), 107--118.

\bibitem{AGM}
Aleksieva Y., Ganchev G., Milousheva V., \textit{On the theory of Lorentz  surfaces with parallel normalized mean curvature vector field  in pseudo-Euclidean 4-space}. J. Korean Math. Soc., \textbf{53}, (2016).

\bibitem{BeCaDu}
Bekta\c{s} B., Canfes E., Dursun U., \emph{On rotational surfaces in pseudo-Euclidean space $\E^4_t$ with pointwise 1-type Gauss map}.
ArXiv 1508.03294v1.

\bibitem{Chen-MM}
Chen B.-Y., \textit{Surfaces with parallel normalized mean curvature vector}.
Monatsh. Math., \textbf{90}, no. 3   (1980), 185--194.

\bibitem{Chen}
Chen B.-Y., \textit{Pseudo-Riemannian geometry, $\delta$-invariants and
applications}. World Scientific Publishing Co. Pte. Ltd.,
Hackensack, NJ, 2011.

\bibitem{Dur}
Dursun U., \emph{On spacelike  rotational surfaces  with pointwise 1-type Gauss map}.
Bull. Korean Math. Soc., \textbf{52}, no. 1  (2015), 301--312.

\bibitem{Dur&Tur1}
Dursun U., Turgay N., \emph{General rotational surfaces in Euclidean space $\E^4$ with pointwise 1-type Gauss map}.
Math. Commun., \textbf{17}  (2012), 71--81.

\bibitem{Dur&Tur2}
Dursun U., Turgay N., \emph{Minimal and pseudo-umbilical rotational surfaces in Euclidean space $\E^4$}.
Mediterr. J. Math. \textbf{10}, no. 1 (2013), 497--506.

\bibitem{GM2}
Ganchev G.,  Milousheva V., \emph{Invariants of lines on surfaces in $\R^4$}. C. R. Acad. Bulgare Sci., {\bf 63}, no. 6 (2010),
835-842.

\bibitem{GM5}
Ganchev G., Milousheva V., \textit{Quasi-minimal rotational surfaces in pseudo-Euclidean
four-dimensional space}. Cent. Eur. J. Math., \textbf{12}, no. 10 (2014),  1586--1601.

\bibitem{GM-TJM}
Ganchev G., Milousheva V., \textit{General rotational surfaces in the 4-dimensional Minkowski space}. Turk. J. Math., 
\textbf{38}, no. 5 (2014), 883--895.

\bibitem{Liu-Liu}
HuiLi L., GuiLi L., \textit{Rotation surfaces with constant mean curvature in 4-dimensional pseudo-Euclidean space}. Kyushu J. Math., \textbf{48}, no. 1  (1994), 35--42. 

\bibitem{Mil}
Milousheva V., \emph{General rotational surfaces in $\R^4$ with
meridians lying in two-dimensional  planes}.
 C. R. Acad. Bulgare Sci., \textbf{63}, no. 3 (2010),
339--348.

\bibitem{Moore1}
Moore C., \emph{Surfaces of rotation in a space of four dimensions}.
Ann. of Math., 2nd Ser., \textbf{21},  no. 2 (1919), 81-93.

\bibitem{Moore2}
Moore C.,   \emph{Rotation surfaces of constant curvature in space of four dimensions}. Bull. Amer. Math. Soc.,
\textbf{26}, no. 10 (1920),  454-460.

\bibitem{Rosca}
Rosca R., \textit{On null hypersurfaces of a Lorentzian manifold}. Tensor (N.S.) \textbf{23}  (1972), 66--74.


\end{thebibliography}
\end{document}